\documentclass[10pt]{amsart}
\usepackage{amsmath,amscd,amssymb,pb-diagram,lamsarrow,pb-lams}
\title{{\large $\boldsymbol{\pi_*}$}-KERNELS OF LIE GROUPS}
\author{Ken-ichi Maruyama}
\address{Department of Mathematics, Faculty of Education, Chiba University, Yayoicho, Chiba, Japan}
\email{maruyama@faculty.chiba-u.jp} 
\begin{document}
\maketitle
\begin{abstract}
We study the group of homotopy classes of self maps of compact Lie groups which induce the trivial homomorphism on homotopy groups. We completely determine the groups for $SU(3)$ and $Sp(2)$. We investigate these groups for simple Lie groups in the cases when Lie groups are $p$-regular or quasi $p$-regular. We apply our results to the groups of self homotopy equivalences.   
\end{abstract}
\newtheorem{prop}{Proposition}[section]
\newtheorem{lemma}[prop]{Lemma}
\newtheorem{theorem}[prop]{Theorem}
\newtheorem{corollary}[prop]{Corollary}
\newtheorem{example}[prop]{Example}
\newtheorem{conjecture}[prop]{Conjecture}
\newtheorem{definition}[prop]{Definition}
\numberwithin{equation}{section}
\begin{center}
{\sc Introduction}
\end{center}
Let $[X, Y]$ be the set of based homotopy classes of maps from a space $X$ to a space $Y$. We denote by $\mathcal Z^n(X, Y)$ the subset of $[X, Y]$ consisting of all homotopy classes which induce the trivial homomorphism on homotopy groups in dimensions less than or equal to $n$. We denote by $\mathcal Z^\infty(X, Y)$ the subset of $[X, Y]$ consisting of all homotopy classes which induce the trivial homomorphism on all homotopy groups.  We write  $\mathcal Z^n(X)$ and $\mathcal Z^{\infty}(X)$ if $X = Y$ for short. In stable theory the set $\mathcal Z^\infty({\bf X}, {\bf Y})$ has been previously considered by Christensen \cite{ch}, where ${\bf X}, {\bf Y}$ are spectra. He calls the elements of $\mathcal Z^\infty({\bf X}, {\bf Y})$ ghosts. Indeed there is a conjecture by Freyd \cite{fr} which states that $\mathcal Z^\infty({\bf X}, {\bf Y})$ is trivial for finite spectra. On the other hand, in the unstable case the situation is quite different. $\mathcal Z^\infty(X, Y)$ is often nontrivial, even infinite for some spaces. 

Let us consider the case where $Y$ is an $H$-space.        
 If $Y$ is an $H$-complex, $[X, Y]$ is an algebraic loop which is actually a group if $Y$ is homotopy associative by a result of James \cite{ja}.  In this case, $\mathcal Z^n(X, Y)$ is a normal subgroup of $[X, Y]$ for any $n$. It is not commutative, but nilpotent in many cases. A main object of this paper is a study of the groups ${\mathcal Z}^n(X)$ for simple Lie groups. We will show that these groups have rather simple structure and can be computable in many cases. We obtain the following theorem.\\
\ \ \\
{\bf Theorem 3.3.}
\begin{gather}
\tag{1} {\mathcal Z}^\infty(SU(3)) = {\mathcal Z}^n(SU(3)) \cong {\bf Z}_{12}~~~~\text{for}~~~~ n \geq 5.\\
\tag{2} {\mathcal Z}^\infty(Sp(2)) = {\mathcal Z}^n(Sp(2)) \cong {\bf Z}_{120}~~~~\text{for}~~~~n \geq 7.
\end{gather}

When $G$ is a compact connected Lie group it is known \cite{ma} that ${\mathcal Z}^\infty(G)$ is equal to ${\mathcal Z}^n(G)$ for some $n$. The smallest such $n$ is written $sz(G)$ and called the stability of  a descending series $\{\mathcal Z^n(G)\}$. We denote  by $lz(G)$ the length of $\{\mathcal Z^n(G)\}$. We also define invariants of $G$ by $sz_p(G) = sz(G_p)$ and $lz_p(G) = lz(G_p)$ respectively. We will show that if $G$ is $p$-regular or  quasi $p$-regular, then $sz_p(G)$  and  $lz_p(G)$ are determined by its rational type:\\
\ \ \\
{\bf Theorem 5.1.} \it  
Let $G$ be a compact connected, simply connected simple Lie group,  $H^*(G;{\bf Q}) \cong \Lambda(x_1,...,x_r)$ with ${\rm deg}x_i = n_i$, $n_1 \leq \cdots \leq n_r$. If $G$ is quasi $p$-regular, then  $sz_p(G) = n_r$ for all $G$ and
\begin{enumerate}
\item[(1)] If $G$ is not isomorphic to $Spin(4n)$
then  $lz_p(G) = r = {\rm rank}~~G$.
\item[(2)] If $G$ is isomorphic to $Spin(4n)$, then $lz_p(G)= r - 1 = {\rm rank}~~G - 1.$ 
\end{enumerate}
\rm   
Here we note that $sz(G) \geq sz_p(G)$ and $lz(G) \geq lz_p(G)$ for an arbitrary prime number $p$ (see  section 4).

 We now briefly summarize the contents of this paper. In section 1 we collect some facts about homotopy sets of $H$-spaces.  In section 2 we determine two numerical invariants $sz(X)$ and $lz(X)$ for some elementary spaces. In section 3 we consider $\mathcal Z^n(G)$ for rank 2 Lie groups $SU(3)$, $Sp(2)$ and $G_2$. We completely determine the group structures of $\mathcal Z^n(G)$ for $SU(3)$ and $Sp(2)$ for all $n$. In section 4 we give $p$-local consideration  to $\mathcal Z^n(G)$ and $\mathcal Z^\infty(G)$ for simple Lie groups. Then we obtain $p$-local information on the invariants $sz(G)$ and $lz(G)$ when $G$ is $p$-regular. In section 5 a result in section 4 is generalized by using quasi $p$-regularity of Lie groups.  In section 6 we give an application of the results in the previous sections. We will study the groups of self homotopy equivalences and their subgroups. The groups are closely related to the groups ${\mathcal Z}^{n}(G)$ when $G$ is a topological group. In particular, we will determine ${\mathcal E}_\#^{\infty}(G)$, the group of homotopy classes which induce the identity on all homotopy groups for $SU(3)$ and $Sp(2)$.  
\section{Preliminaries}
In this section we fix our notation and give some general results. Throughout this paper all spaces are connected and based. All maps and homotopies preserve base points.  We do not distinguish notationally between a map and its homotopy class. For spaces $X$ and $Y$, we let $[X, Y]$ denote the set of homotopy classes from $X$ to $Y$. As we noted in our introduction, if $Y$ is a homotopy associative $H$-space then $[X, Y]$ is a group with the binary operation obtained from the $H$-structure of $Y$ \cite{ja}. We denote this operation ``+''. That is,
$$ f + g = \mu_Y\circ (f \times g) \circ \Delta,$$ 
where $\mu_Y : Y \times Y \to Y$ is the $H$-structure of $Y$, $f, g \in [X, Y]$,  and $\Delta : X \to X \times X$ is the diagonal map. 
Then for any map $h : W \to X$,
\begin{equation}
(f + g) \circ h = f\circ h + g\circ h\hspace{0.3cm}in\hspace{0.2cm}[W, Y]. 
\end{equation}
Let ${\mathcal Z}^n(X, Y) \subset [X, Y]$ denote the subset which consists of elements $f : X \to Y$ such that $f_* = 0 : \pi_i(X) \to \pi_i(Y)$ for $i \leq n$. We denote by ${\mathcal Z}^\infty(X, Y) \subset [X, Y]$ the subset which consists of elements $f : X \to Y$ such that $f_* = 0 : \pi_i(X) \to \pi_i(Y)$ for all $i$. We write ${\mathcal Z}^n(X)$ for ${\mathcal Z}^n(X, X)$ and ${\mathcal Z}^\infty(X)$ for ${\mathcal Z}^\infty(X, X)$. We call ${\mathcal Z}^n(X)$ {\it the nth $\pi_*$-kernel of $X$}.
\begin{prop}
Let $Y$ be a homotopy associative $H$-space. Then ${\mathcal Z}^n(X, Y)$ is a normal subgroup of the group $[X, Y]$.
\end{prop}
\begin{proof}
By (1.1), we obtain that $(f + g)_*(x) = f_*(x) + g_*(x)$ for $f, g \in [X, Y]$ and $x \in \pi_*(X)$. Thus the result follows.   
\end{proof}
\begin{corollary}
If $X$ is a connected finite homotopy associative $H$-space, then ${\mathcal Z}^n(X)$ is a nilpotent group. 
\end{corollary}
\begin{proof}
It is known that $[W, X]$ is a nilpotent group if $W$ is finite-dimensional (cf. \cite{wh}). The result follows from Proposition 1.1.
\end{proof}  
%
%
\section{numerical invariants}
In this section we define two numerical invariants related to ${\mathcal Z}^*(X)$. We will investigate these invariants in later sections. 
\begin{definition}
If there exists an integer $t \geq 1$ such that ${\mathcal Z}^n(X) = {\mathcal Z}^{n+1}(X)$ for $n \geq t$, the smallest such $t$ is written $sz(X)$. If no such integer exists we write $sz(X) = \infty$. We call $sz(X)$ the {\it stability} of ${\mathcal Z}^*(X)$.
\end{definition}
\begin{definition} 
The number of strict inclusions in $[X, X] \supset {\mathcal Z}^1(X) \supset {\mathcal Z}^2(X) \supset \cdots$ is denoted by $lz(X)$. We call $lz(X)$ the length of ${\mathcal Z}^*(X)$. 
\end{definition}
\noindent
{\bf Remark}. Clearly $sz(X) \geq lz(X)$. We may consider that $lz(X)$ reflects algebraic properties of ${\mathcal Z}^*(X)$ while $sz(X)$ is a geometric invariant.\\
\ \\
We introduce two lemmas which will be needed in the sequel. 
\begin{lemma}
Let $\prod_{i=1}^r S^{n_i}$ be the product of spheres with $n_1 \leq \cdots \leq n_r$. Then $ sz(\prod_{i=1}^r S^{n_i}) = n_r$, and $lz(\prod_{i=1}^r S^{n_i}) = \#\{n_1,...,n_r\}$ $($= the number of the distinct integers in $\{n_1,...,n_r\})$.
\end{lemma}
\begin{proof}
The inclusion $\vee_{i=1}^r S^{n_i} \to \prod_{i=1}^r S^{n_i}$ induces an epimorphism on homotopy groups. It follows that $f \in {\mathcal Z}^{n_r}(\prod_{i=1}^r S^{n_i})$ implies $f\mid_{ \vee_{i=1}^r S^{n_i}} \simeq *$ and thus $f \in  {\mathcal Z}^\infty(\prod_{i=1}^r S^{n_i})$. Therefore we obtain the first equality. For the second equality, note that there exists a map $ * \times id : \prod_{i=1}^j S^{n_i} \times \prod_{i=j+1}^r S^{n_i} \to \prod_{i=1}^j S^{n_i} \times \prod_{i=j+1}^r S^{n_i}$ which belongs to ${\mathcal Z}^{n_j}(\prod_{i=1}^r S^{n_i})$ but does not belong to ${\mathcal Z}^{n_{j+1}}(\prod_{i=1}^r S^{n_i})$ if $n_j \not= n_{j+1}$. Thus $lz(\prod_{i=1}^r S^{n_i}) \geq \#\{n_1,...,n_r\}$. Next assume that there exists $n$ such that
$${\mathcal Z}^{n_j}(\prod_{i=1}^r S^{n_i}) \supsetneqq {\mathcal Z}^{n}(\prod_{i=1}^r S^{n_i}) \supsetneqq {\mathcal Z}^{n_{j+1}}(\prod_{i=1}^r S^{n_i})$$
for some $ 1 \leq j \leq r$. Then there exists $f \in {\mathcal Z}^{n_j}(\prod_{i=1}^r S^
{n_i})$ such that $f_* \not= 0$ on $\pi_n(\prod_{i=1}^r S^{n_i})$. Hence the composition 
$p_{k'}fi_k : S^{n_k} \to S^{n_{k'}}$ is nontrivial for some $k$, $k'$ $\leq j$, where $i_k : S^{n_k} \to \prod_{i=1}^r S^{n_i}$ and $p_{k'} : \prod_{i=1}^r S^{n_i} \to S^{n_{k'}}$ are the inclusion and the projection. This is a contradiction, because the restriction $f$ to $\vee_{i=1}^j S^{n_i}$ is null homotopic. Therefore we obtain $lz(\prod_{i=1}^r S^{n_i}) \leq \#\{n_1,...,n_r\}$ and the result follows.    
  
\end{proof}
We obtain a variation of the above lemma as follows.
\begin{lemma} Let $\prod_{i=1}^r B_i$ be a product space such that ${\mathcal Z}^\infty(B_i, B_j)$ $=$ $[B_i, B_j]$, for $i \not= j$.
\begin{enumerate}
\item[(1)] Assume that $sz(B_i) = s_i$. Then $sz(\prod_{i=1}^r B_i)$ $=$ $\max\{s_i\}$.
\item[(2)] Assume that $ lz(B_i) = l_i$ and ${\mathcal Z}^{i_k}(B_i) \not= {\mathcal Z}^{i_k + 1}(B_i)$ for $k = 1,...,l_i$. Then $lz(\prod_{i=1}^r B_i)$ $=$ $\#\{i_k\}$, $i = 1,...,r$, $k = 1,...,l_i$.
\end{enumerate}
\end{lemma}
\begin{proof}
The proof is similar to that of Lemma 2.3.
\end{proof}
\section{Lie groups of rank 2} 
In this section we study ${\mathcal Z}^n(X)$ or ${\mathcal Z}^\infty(X)$ for $SU(3)$, $Sp(2)$ and $G_2$. The self homotopy sets $[X, X]$ have been determined by Mimura and \=Oshima \cite{mio} when $X = SU(3), Sp(2)$ and $[G_2, G_2]$ has been determined by \=Oshima up to group extension in \cite{os}. They have obtained the following results.
\begin{theorem}[{\rm \cite{mio} \cite{os}}] Relative to the standard multiplications, 
\begin{align*}
[SU(3), SU(3)] & \cong \Psi(12, 1)\\
[Sp(2), Sp(2)] & \cong \Psi(120, 12)
\end{align*}
The following sequence is exact.
\begin{equation}
 0  \to \pi_{14}(G_2) \xrightarrow{q_{14}^*}[G_2, G_2] \to \Psi(2,1) \to 0 
\end{equation}
Here $\Psi(m, n)$ is the group with generators $x, y, z$ and relations 
$$ xz = zx,~ yz = zy,~ z^m=1,~ [x, y] = z^n.$$  
\end{theorem}
By Theorem 3.1 and the results of \cite{mio}we obtain
\begin{theorem}
\begin{gather}
\tag{1} {\mathcal Z}^n(SU(3)) \cong
\begin{cases}
[SU(3), SU(3)]& \quad \text{for $n < 3$}\\ 
{\bf Z} \oplus {\bf Z}_{12},& \quad \text{for $3 \leq n < 5$}\\
{\bf Z}_{12}, & \quad \text{for $n = 5$}
\end{cases}\\
{\mathcal Z}^n(Sp(2)) \cong 
\begin{cases}
[Sp(2), Sp(2)]& \quad \text{for $n < 3$}\\
\tag{2} {\bf Z} \oplus {\bf Z}_{120},& \quad \text{for $3 \leq n < 7$}\\
{\bf Z}_{120}, & \quad \text{for $n = 7$}
\end{cases}
\end{gather}
The torsion subgroup is generated by $z$ in the both cases.
\end{theorem}
\begin{proof}
We have $\pi_5(SU(3)) \cong {\bf Z}$, $\pi_8(SU(3)) \cong {\bf Z}_{12}$ and $\pi_{7}(Sp(2)) \cong {\bf Z}$, $\pi_{10}(Sp(2)) \cong {\bf Z}_{120}$ (see \cite{mito}). By \cite{mio}, $[SU(3), SU(3)]$ and $[Sp(2), Sp(2)]$ are generated by
\begin{align}
x &: X \xrightarrow{p} S^n \xrightarrow{a}  X\\
y &= id : X \to X\\
z &: X \xrightarrow{q}  S^{n + 3} \xrightarrow{b}  X,
\end{align} 
where $n = 5$ for $X = SU(3)$ and $n = 7$ for $X = Sp(2)$, $p : X \to S^n$ is the bundle projection, $q : X \to S^{n+3}$ is the projection to the top cell,  $a$ and $b$ are generators of the homotopy groups. 
 Let us consider the $SU(3)$ case. By \cite{mi} $\pi_5(SU(3))$ is generated by the element $[2\iota_5]$ such that $p_*([2\iota_5]) = 2\iota_5$, where $\iota_5$ is the identity class of  $\pi_5(S^5)$. Thus $a$ in (3.2) is equal to $[2\iota_5]$. Since we have 
$$ x_*([2\iota_5]) = [2\iota_5] \circ p \circ [2\iota_5] = 2[2\iota_5], $$
$x$ induces a nontrivial homomorphism on $\pi_5(SU(3))$.  Clearly $x$ induces the trivial homomorphisms on $\pi_3(SU(3))$ and $\pi_4(SU(3))$, and hence $x$ $\in {\mathcal Z}^4(SU(3)$ and $x \notin {\mathcal Z}^5(SU(3))$. The map $z$ induces the zero homomorphism on homotopy groups $\pi_i(SU(3))$, $i \leq 8$ as $z$ factors through $S^8$. Therefore, ${\mathcal Z}^n(SU(3))$ for $n < 5$ are generated by $x$ and $z$, ${\mathcal Z}^n(SU(3))$ for $5 \leq n \leq 8$ are generated by $z$. The group structures are obtained by Theorem 3.1. 

The case where $X = Sp(2)$ is proved by the same methods.      
\end{proof}
Next we consider ${\mathcal Z}^\infty(G)$ for $G = SU(3)$ and $Sp(2)$. We can completely determine these groups. Our result is as follows.
%
\begin{theorem} 
\begin{gather}
\tag{1} {\mathcal Z}^\infty(SU(3)) = {\mathcal Z}^n(SU(3)) \cong {\bf Z}_{12}~~~~\text{for}~~~~ n \geq 5.\\
\tag{2} {\mathcal Z}^\infty(Sp(2)) = {\mathcal Z}^n(Sp(2)) \cong {\bf Z}_{120}~~~~\text{for}~~~~n \geq 7.
\end{gather}
\end{theorem}
\begin{proof}
By Theorem 3.2, ${\mathcal Z}^5(SU(3))$ is generated by the commutator $[x, y]$. Since the group structure of $[SU(3), SU(3)]$ is induced from the multiplication of $SU(3)$, $z = [x, y]$ induces $z_* = x_* + y_* - x_* - y_* = 0$ on homotopy groups, and hence we obtain (1). Next we prove (2). Let $\omega_n$ denote the generator of the group $\pi_{n+3}(S^n) \cong {\bf Z}_{24}$, $n \geq 5$. According to \cite{mi} Proposition 4.1, $Sp(2) \wedge Sp(2)/S^6$ has the following homotopy type.
$$ S^{10} \vee S^{10} \cup_\xi C(S^{13} \vee S^{19}) \vee S^{13} \cup_{2\omega_{13}} e^{17} \vee S^{13} \cup_{2\omega_{13}} e^{17},$$
where $\xi = 2\omega_{10} + 4\omega_{10} + \delta$, and $\delta$ is a Whitehead product. From \cite{mito} we see that $[\nu_7]\circ 4\omega_{10} = 0$. 
 Thus the generator $[\nu_7]$ of $\pi_{10}(Sp(2)) \cong {\bf Z}_{120}$ extends to $Sp(2) \wedge Sp(2)/S^6$. We denote by $\overline{[\nu_7]}$ an extension of $[\nu_7]$:
\begin{eqnarray*}
\begin{diagram}
\node{Sp(2) \wedge Sp(2)/S^6} \arrow{se,t}{\overline{[\nu_7]}}\\
\node{S^{10} \vee S^{10}}\arrow{n} \arrow{e,b}{0 \vee [\nu_7]} \node{Sp(2)}
\end{diagram}
\end{eqnarray*}
Now we consider the composition:
\begin{equation}
 Sp(2) \xrightarrow{\bar{\Delta}} Sp(2) \wedge Sp(2) \to Sp(2) \wedge Sp(2)/S^6 \xrightarrow{\overline{[\nu_7]}} Sp(2),
\end{equation}
where $\bar{\Delta}$ is the reduced diagonal map. By connectivity  $\bar{\Delta}$ induces a map $k$ which makes the following diagram commutative.
\begin{eqnarray*}
\begin{diagram}
\node{Sp(2)} \arrow{e,t}{\bar{\Delta}} \arrow{s,r}{q}
\node{Sp(2) \wedge Sp(2)} \arrow{e} \node{Sp(2) \wedge Sp(2)/S^6} \arrow{e,t}{\overline{[\nu_7]}} \node{Sp(2)}\\
\node{S^{10}} \arrow{nee,r}{k} 
\end{diagram}
\end{eqnarray*}
The map $k$ is homotopic to $i_1 + i_2 : S^{10} \to Sp(2) \wedge Sp(2)/S^6$, where $i_1$ and $i_2$ are the inclusions. As $\overline{[\nu_7]}k = [\nu_7]$, the composition (3.5) is homotopic to the generator $z$ of ${\mathcal Z}^7(Sp(2)) \cong {\bf Z}_{120}$. Since $\bar{\Delta} \in {\mathcal Z}^\infty(Sp(2), Sp(2) \wedge Sp(2))$, we obtain $z \in {\mathcal Z}^\infty(Sp(2))$. 
\end{proof}

\begin{corollary}
$sz(SU(3)) = 5$, $lz(SU(3)) = 2$. $sz(Sp(2)) = 7$, $lz(Sp(2)) = 2$.
\end{corollary}
We obtain a corresponding result to Theorem 3.2 for $G_2$  as follows. 
\begin{theorem}
\begin{equation*}
{\mathcal Z}^n(G_2) \cong
 \begin{cases}
[G_2, G_2], & \quad \text{for $n < 3$}\\
{\bf Z}_2 \oplus {\bf Z}_2 \oplus {\bf Z}_8 \oplus {\bf Z}_{21} \oplus {\bf Z},&  \quad \text{for $3 \leq n < 11$}\\
{\bf Z}_2 \oplus {\bf Z}_2 \oplus {\bf Z}_8 \oplus {\bf Z}_{21}, & \quad \text{for $11 \leq n \leq 14$}
\end{cases}
\end{equation*}
\end{theorem}
\begin{proof}
In the exact sequence (3.1) in Theorem 3.1, $\pi_{14}(G_2)$ is isomorphic to ${\bf Z}_8 \oplus {\bf Z}_2 \oplus {\bf Z}_{21}$, and the group $\Psi(2, 1) \cong [G_2^{11}, G_2] $ is generated by the elements as follows \cite{os}:
\begin{align}
 i_{11} & : G_2^{11} \to G_2\\
 q_{11,6}^*\gamma'& : G_2^{11} \xrightarrow{q_{11,6}} G_2^{11}/G_2^{6} \to G_2 \\
 q_{11}^*j_*[\nu_5^2] & : G_2^{11} \xrightarrow{q_{11}} S^{11} \xrightarrow{j} G_2
\end{align}   
There exists a short exact sequence
\begin{equation}
 0 \to \pi_{14}(G_2) \xrightarrow{q_{14}^*} {\mathcal Z}^n(G_2) \xrightarrow{i_{11}^*} {\bf Z} \oplus {\bf Z}_2 \to 0 
\end{equation}
for $ 3 \leq n \leq 10$.
If $n \leq 7$, then the right term is generated by $q_{11,6}^*\gamma'$ and $q_{11}^*j_*[\nu_5^2]$, and if $ 7 < n \leq 14$ then it is generated by  $q_{11}^*j_*[\nu_5^2]$ and $\epsilon q_{11,6}^*\gamma'$, $\epsilon$ = 1 or 2. This can be seen as follows.
 Since the inclusion map induces an epimorphism $\pi_i(G^{11}) \to \pi_i(G)$ for $i \leq 13$, any extensions of $q_{11,6}^*\gamma'$ and $q_{11}^*j_*[\nu_5^2]$ induce the trivial homomorphism on $\pi_i(G_2)$ for $i \leq 7$. 
  Let $\alpha$ be an extension of $q_{11,6}^*\gamma'$, that is, $i_{11}^*\alpha = q_{11,6}^*\gamma'$. Then by \cite{os} Theorem 2.2, $i_{11}^*[1, \alpha] = [i_{11}, q_{11,6}^*\gamma'] = q_{11}^*j_*[\nu_5^2]$. Thus $[1, \alpha]$ is an extension of $q_{11}^*j_*[\nu_5^2]$ and induces the trivial homomorphism on all homotopy groups as we saw in the proof of Theorem 3.3. Moreover it is easy to see that $2\alpha$ induces the trivial homomorphism on $\pi_8(G_2)$ and $\pi_9(G_2)$. Hence we obtain the exact sequence (3.9) (note that $\pi_{10}(G_2) = 0$).  Now by \cite{os}, 
  $$[\alpha, [1,\alpha]] = 0$$
Therefore, ${\mathcal Z}^n(G_2)$, $3 \leq n \leq 10$ are abelian groups. On the other hand, 
$$ G_2 \to G_2/G_2^{9} \simeq S^{11} \vee S^{14} \to S^{11} \xrightarrow{j_*[\nu_5^2]} G_2 $$is also an extension of $q_{11}^*j_*[\nu_5^2]$ of order 2.  Therefore, the above exact sequence (3.9) splits.  
 By \cite{oka} Lemma 5.8  $\bar{q}^*_{11,6}\gamma = 4\gamma'$, where $\gamma : S^{11} \to G_2$ is a generator of the direct summand ${\bf Z}$ of $\pi_{11}(G_2)$, and $\bar{q}_{11,6} : G_2^{11}/G_2^6 \to S^{11}$ is the projection map. Therefore an arbitrary extension of  $q_{11,6}^*\gamma'$ to  ${\mathcal Z}^n(G_2)$ induces a nontrivial homomorphism on $\pi_{11}(G_2)$. The homotopy groups $\pi_{n}(G_2)$ are zero for $n = 12, 13$. Now by definition we have $q_{14}^*(\pi_{14}(G_2)) \subset {\mathcal Z}^{14}(G_2)$. As $[1, \alpha]$ induces the trivial homomorphism on all the homotopy groups, the result follows.   
\end{proof}
\section{$p$-local cases}
For a set of prime numbers $P$ and a nilpotent space (or group) $X$, we denote by $X_P$ its $P$-localization. In particular, let $X_0$ denote the rationalization of a nilpotent space (or group) $X$. If $X$ is a finite $H$-space, its rational cohomology $H^*(X; {\bf Q})$ is isomorphic to $\Lambda(x_1,...,x_\ell)$, the exterior algebra over ${\bf Q}$ with deg$x_i$ odd. We call $\ell$ the rank of $X$ and (deg$x_1$,...deg$x_\ell$) the type of $X$. It is known that $X_0$ is homotopy equivalent to  $\prod_{i=1}^{\ell}S^{{\rm deg}x_i}_0$.
Let us begin with definitions of the local stability and the local length of ${\mathcal Z}^*(X)$. 
\begin{definition}
 Let $X$ be a nilpotent space, $P$ a set of prime numbers . We define two invariants of $X$ by  
$sz_P(X) = sz(X_P)$ and $lz_P(X) = lz(X_P)$. 
\end{definition} 
%
We recall two results of \cite{ma}.
\begin{lemma}[\cite{ma}]Let $X$ be a homotopy associative finite $H$-space, $P$ a set of prime numbers.
\begin{enumerate}
\item[(1)] $sz(X)$ and $lz(X)$ are finite.
\item[(2)] ${\mathcal Z}^n(X)_P \cong {\mathcal Z}^n(X_P)$ for any $n >0$ and $n = \infty$. 
\end{enumerate}
$($Thus $sz_P(X)$ and $lz_P(X)$ are finite.$)$
\end{lemma}
\noindent
{\bf Remark}. We remark that if $X$ is a homotopy associative finite $H$-space, then $sz(X) = \max_p\{sz_p(X)\}$ and $lz(X) = \max_p\{lz_p(X)\}$\vspace{0.3cm}.\par 

We now completely determine the odd primary part of the group ${\mathcal Z}^*(G_2)$.
\begin{theorem}
${\mathcal Z}^\infty(G_2)_{(odd)} = {\mathcal Z}^{n}(G_2)_{(odd)} \cong {\bf Z}_{21}$ for $n \geq 11$.
\end{theorem}
\begin{proof}
Since ${\mathcal Z}^{11}(G_2)_{(odd)} \cong {\bf Z}_{21}$, we only have to consider the 3 and 7 primary cases.  First we consider the case when $p =3$. It is known that $(G_2)_3$ is equivalent to the principal $S^3$-bundle over $S^{11}$ with the characteristic element $\alpha_2$. Thus $G_2$ at $3$ has the following cell structure:
$$ S^3 \cup_{\alpha_2} e^{11} \cup_\rho S^{14} $$
By results of James \cite{ja2} $\Sigma \rho$ is homotopic to $\Sigma i \circ J(\chi)\), where $i : S^3 \to G_2$, and $J(\chi)$ is the $J$-image of the characteristic element $\chi \in \pi_{10}(O(4))$ of the bundle. Therefore $J(\chi) = \alpha_1(4) \circ \alpha_2(7)$ or 0 by \cite{to}. By using the fact that $\alpha_1(7) \circ \alpha_2(10) = 0$ \cite{to}, $\Sigma^4 \rho $ is trivial. Consider the following homotopy commutative diagram.  
\begin{eqnarray}
\begin{diagram}
\node{G_2} \arrow{e,t}{(1 \times p_{11})\Delta} \arrow{se,r}{((1 \times p_{11})\Delta)'}
\node{G_2 \times S^{11}} \arrow{e,t}{\wedge} \node{G_2 \wedge S^{11}}\\
\node{} \node{(G_2 \times S^{11})^{14}} \arrow{n} \arrow{e,t}{\wedge'} \node{S^3 \wedge S^{11}} \arrow{n}
\end{diagram}
\end{eqnarray}   
where $\Delta$ is the diagonal map, $\wedge$ is the projection and $p_{11}$ is the bundle projection, $((1 \times p_{11})\Delta)'$ and $\wedge'$ are maps induced from $(1 \times p_{11})\Delta$ and $\wedge$ by dimensional reasons. Clearly $\wedge'\circ ((1 \times p_{11})\Delta)'$ is homotopic to $q$, the projection map to the top cell. As we saw that $\Sigma^4\rho = 0$, $(G_2 \wedge S^{11})_3 \simeq (S^{14} \cup_{\alpha_2} e^{22} \vee S^{25})_3$. From the result in the previous section we see that ${\mathcal Z}^{11}(G_2)_3 \cong {\bf Z}_3$ is generated by the composition:
$$ q^*(\alpha_3) : G_2 \to S^{14} \xrightarrow{\alpha_3} G_2 $$ 
Since $\pi_{21}(G_2) = 0$ \cite{mi2}, $\alpha_3$ extends to $S^{14} \cup_{\alpha_2} e^{22}$, namely we obtain the following commutative diagram localized at 3.
\begin{eqnarray*}
\begin{diagram}
\node{S^{14} \cup_{\alpha_2} e^{22}} \arrow{se,t}{\overline{\alpha_3}}\\
\node{S^{14}}\arrow{n} \arrow{e,b}{\alpha_3} \node{G_2}\\
\node{S^{21}}\arrow{n,l}{\alpha_2}
\end{diagram}
\end{eqnarray*}
We have the following commutative diagram at 3:
\begin{equation}
\begin{diagram}
\node{G_2} \arrow{e,t}{(1 \times p_{11})\Delta} \arrow{se,r}{((1 \times p_{11})\Delta)'}  
\node{G_2 \times S^{11}} \arrow{e,t}{\wedge} \node{G_2 \wedge S^{11}} \arrow{se,t}{\overline{\alpha_3} \vee 0}\\
\node{} \node{(G_2 \times S^{11})^{14}} \arrow{n} \arrow{e,t}{\wedge'} \node{S^3 \wedge S^{11}} \arrow{n} \arrow{e,t}{\alpha_3} \node{G_2}
\end{diagram}
\end{equation}
Here we should note that $G_2 \wedge S^{11}$ is 3-equivalent to $S^{14} \cup_{\alpha_2} e^{22} \vee S^{25}$.
In (4.2), the lower composition is homotopic to $q^*(\alpha_3)$, and the upper composition is an element of ${\mathcal Z}^\infty(G_2)$ since $\wedge$ induces the trivial map on homotopy groups. Thus we obtain the result for $p = 3$. For $p = 7$, $G_2$ is 7-equivalent to the product $S^3 \times S^7$, and hence it easily follows that ${\mathcal Z}^\infty(G_2)_7 \cong {\mathcal Z}^\infty((G_2)_7) \cong {\mathcal Z}^\infty((S^3 \times S^{11})_7) \cong {\mathcal Z}^\infty(S^3 \times S^{11})_7 \cong {\mathcal Z}^{11}(S^3 \times S^{11})_7 \cong {\bf Z}_7$ by Lemma 2.3 and Lemma 4.2.  
\end{proof}
\begin{corollary}
$sz_{(odd)}(G_2) = 11$, $lz_{(odd)}(G_2) = 2$ 
\end{corollary}
Next we will study ${\mathcal Z}^\infty(G)_p$, $sz_p(G)$ and $lz_p(G)$ for compact Lie groups. By using Lemma 4.2 we can obtain a finiteness property of ${\mathcal Z}^\infty(G)$ for compact Lie groups $G$.
\begin{prop}
Let $G$ be a compact connected, simply connected simple Lie group, then ${\mathcal Z}^\infty(G)$ is finite if and only if $G$ is isomorphic to one of the following groups:
\begin{enumerate}
\item[(1)] $SU(n)$ $n < 8$
\item[(2)] $Spin(n)$ $n < 29$
\item[(3)] $Sp(n)$ $n < 14$
\item[(4)] Exceptional Lie groups other than $E_6$
\end{enumerate} 
\end{prop}
\begin{proof}
As $G_0$ is homotopy equivalent to  $\prod_{i=1}^r S^{n_i}_0$, where $n_1 \leq \cdots \leq n_r$ and $n_i$ is an odd integer, ${\mathcal Z}^\infty(G)_0$ is isomorphic to ${\mathcal Z}^\infty(\prod_{i=1}^r S^{n_i}_0)$. The latter group is isomorphic to the subgroup of $[\prod_{i=1}^r S^{n_i}_0$, $\prod_{i=1}^r S^{n_i}_0]$ $\cong$ $\oplus_{i=1}^r H^{n_i}(\prod_{i=1}^r S^{n_i}_0)$ generated by decomposable elements of dimensions $\leq n_r$. Our result is derived from the rational types and observation of the Poincar\'e polynomials of these groups. 
\end{proof}
{\bf Remark}. 
\begin{enumerate}
\item[(1)] If $G$ is a compact connected, simple Lie group then ${\mathcal Z}^\infty(G)$ is finite if and only if ${\mathcal Z}^\infty(\tilde{G})$ is finite, where $\tilde{G}$ is the universal cover of $G$. For example, ${\mathcal Z}^\infty(SO(n))$ is finite if and only if ${\mathcal Z}^\infty(Spin(n))$ is finite.
\item[(2)] We also see easily that ${\mathcal Z}^\infty(U(n))$ is finite if and only if $n < 5$.
\end{enumerate}
\ \ \\
In the above proposition we used $0$-regularity of Lie groups.  Next we apply $p$-regularity of Lie groups to obtain information on $sz(G)$ and $lz(G)$.      
\begin{prop}
Let $G$ be a compact connected, simply connected simple Lie group. If $G$ is $p$-regular, namely $G_p$  is homotopy equivalent to a product space $\prod_{i=1}^r S^{n_i}_p$, where $n_1 \leq \cdots \leq n_r$ and $n_i$ is an odd integer, then $sz_p(G) = sz_0(G)$, $lz_p(G) = lz_0(G).$ Thus $sz_p(G) = n_r$ for all $G$ and 
\begin{enumerate}
\item[(1)] If $G$ is not isomorphic to $Spin(4n)$
then  $lz_p(G) = r = {\rm rank}~~G$.
\item[(2)] If $G$ is isomorphic to $Spin(4n)$, then $lz_p(G)= r - 1 = {\rm rank}~~G - 1$ 
\end{enumerate}
\end{prop}
\begin{proof}
It is clear from Lemma 2.3 and Lemma 4.2 that $sz_p(G) = sz_0(G) = n_r$ and $lz_p(G) = lz_0(G)$ for any $p$-regular Lie group $G$.   
Let us consider the case $G = Spin(2n)$. If $Spin(2n)$ is $p$-regular, then
$$ Spin(2n)_p \simeq \prod_{i=2}^n S^{4i-5}_p \times S_p^{2n-1}.$$
Therefore
\begin{equation*}
 lz_p(Spin(2n)) =\#\{3, 7,...,4n-5, 2n-1\} =
\begin{cases}
  n& \text{if $n$ is odd},\\
  n-1& \text{if $n$ is even}
\end{cases}
\end{equation*}  
The remaining cases are similar.
\end{proof}
\noindent
{\bf Remark} \begin{enumerate}
\item[(1)]
By the results of \cite{ku} and \cite{se} we obtain the following\par
\begin{enumerate}
\item[]If $p > n-1$, $sz_p(SU(n)) = 2n-1$ and $lz_p(SU(n)) = n-1$ 
\item[]If $p > 2n-1$, $sz_p(Sp(n)) = 4n-1$ and $lz_p(Sp(n)) = n$  
\item[]If $p > n-2$, $sz_p(Spin(2n-1)) = 4n-5$ and $lz_p(Spin(2n-1)) = n-1$
\item[]If $p > n-2$ and $n$ odd, $sz_p(Spin(2n)) = 4n-5$ and $lz_p(Spin(2n)) = n$  
\item[]If $p > n-2$ and $n$ even, $sz_p(Spin(2n)) = 4n-5$ and $lz_p(Spin(2n)) = n -1$ 
\item[]If $p > 5$, $sz_p(G_2) = 11$ and $lz_p(G_2) = 2$ 
\item[]If $p > 11$, $sz_p(F_4) = 23$, $sz_p(E_6) = 23$ and $lz_p(F_4) = 4$, $lz_p(E_6) =6$
\item[]If $p > 17$, $sz_p(E_7) = 35$ and $lz_p(E_7) = 7$ 
\item[]If $p > 29$, $sz_p(E_8) = 59$ and $lz_p(E_8) = 8$ \vspace{0.3cm}
\end{enumerate}
\item[(2)] Since $SO(n)_p \simeq Spin(n)_p$ for an odd prime $p$, we obtain in this case
$$ sz_p(SO(n)) = sz_p(Spin(n)),~~~~lz_p(SO(n)) = lz_p(Spin(n)).$$
\item[(3)] It is known that $U(n)$ is homeomorphic to $S^1 \times SU(n)$. There exists the following exact sequence 
$$ 0 \to [S^1 \wedge SU(n), SU(n)] \to {\mathcal Z}^i(S^1 \times SU(n)) \to {\mathcal Z}^i(SU(n)) \to 0$$  
for $i > 1$. Thus $sz(U(n)) = sz(SU(n))$ for $n > 1$, and $lz(U(n)) = lz(SU(n)) + 1$.
\end{enumerate}
\section{Quasi $p$-regularity}
In the previous section we studied the group ${\mathcal Z}^*(G)_p$ for the case where $G$ is $p$-regular. In this section we will generalize Proposition 4.6. We will show that $p$-regularity can be replaced by quasi $p$-regularity. Quasi regularity of Lie groups has been studied in \cite{oka-1} and \cite{mito}.  Let  $B_n(p)$ be the $S^{2n + 1}$-bundle over $S^{2n+2p-1}$ determined by the generator $\alpha_1$ of $\pi_{2n + 2p -2}(S^{2n+1}) \simeq {\bf Z}_p$. $B_n(p)$  has a cell structure of the form
$$ S^{2n +1} \cup_{\alpha_1} e^{2n + 1 + 2(p-1)} \cup e^{4n + 2 + 2(p-1)}. $$
A Lie group $G$ is called quasi $p$-regular if and only if $G$ is homotopy equivalent to a product of spheres and $B_n(p)$'s at a prime $p$:
$$ \prod S^{m_i}_p \times \prod B_{n_j}(p)_p \simeq G_p $$
\begin{theorem}
Let $G$ be a compact connected, simply connected simple Lie group,  $H^*(G;{\bf Q}) \cong \Lambda(x_1,...,x_r)$ with ${\rm deg}x_i = n_i$, $n_1 \leq \cdots \leq n_r$. If $G$ is quasi $p$-regular, then  $sz_p(G) = n_r$ for all $G$ and
\begin{enumerate}
\item[(1)] If $G$ is not isomorphic to $Spin(4n)$
then  $lz_p(G) = r = {\rm rank}~~G$.
\item[(2)] If $G$ is isomorphic to $Spin(4n)$, then $lz_p(G)= r - 1 = {\rm rank}~~G - 1.$ 
\end{enumerate}
\end{theorem}
In the remaining of this section all spaces are localized at a prime number $p$ which we are considering. For example, $B_n(p)$ means $B_n(p)_p$.  
We quote from \cite{to}
\begin{prop}
Let $p$ be an odd prime.  Then\\
$\pi_{2n + 1 + 2k(p-1) -2}(S^{2n + 1}) \cong {\bf Z}_p $ for $1 \leq n < k$, and $ k = 2,\cdots, p-1$\\
$\pi_{2n + 1 + 2k(p-1) -1}(S^{2n + 1}) \cong {\bf Z}_p$ for $ k = 1, 2,..., p-1,$ $n \geq 1$,\\
$\pi_{2n + 1 + t}(S^{2n + 1}) = 0$ otherwise for $t < 2p(p-1) -2$.     
\end{prop}
\begin{lemma}
Let $p$ be an odd prime. If $[S^{2n + 1} \cup_{\alpha_1} e^{2n + 2p -1}, X] = 0,$ then $q^* : [S^{4n + 2p}, X] \to [B_n(p), X]$ is surjective.
\end{lemma}
\begin{proof}
The result is clear from the following exact sequence.
\begin{equation}
\begin{diagram}
\node{[S^{4n + 2p}, X]} \arrow{e,t}{q^*} \node{[B_n(p), X]} \arrow{e,t}{j^*} \node{[S^{2n + 1} \cup_{\alpha_1} e^{2n + 2p -1}, X]}
\end{diagram}
\end{equation}
Here the sequence is induced by the cofibration.
\end{proof}
\begin{lemma}
Let $p$ be an odd prime , $n$ a positive integer less than $p$ and $m$ be a positive integer such that $n \not= m$ and $n - m + p -1 \not= 0$. Let $q : B_n(p) \to S^{4n + 2p}$ be the projection map to the top cell. Then 
\begin{gather}
\tag{1} q^* : [S^{4n + 2p}, S^{2m+1}] \to [B_n(p), S^{2m+1}].\\ 
\tag{2} q^* : [S^{4n + 2p}, B_m(p)] \to [B_n(p), B_m(p)].
\end{gather}
are surjective.
\end{lemma}
\begin{proof}
First note that $2(n-m) < 2p -3$ and $2n + 1 < 2m + 2p -1$. It follows that
\begin{eqnarray}
\pi_{2n+1}(S^{2m+1}) = 0,~~~~~~
\pi_{2n + 1}(S^{2m + 2p -1}) = 0
\end{eqnarray}
By using the fibration 
\begin{eqnarray}
S^{2m+1} \to B_m(p) \to S^{2m + 2p -1} 
\end{eqnarray}
we also have $\pi_{2n+1}(B_m(p))= 0$, and hence we obtain the exact sequence
\begin{eqnarray}
 [S^{2n +2}, X] \xrightarrow{\alpha_1^*} [S^{2n + 2p -1}, X] \to [S^{2n + 1} \cup_{\alpha_1} e^{2n + 2p -1}, X] \to 0
\end{eqnarray} 
for $X$ = $S^{2m+1}$ or $B_n(p)$.
As $2n - 2m + 2p -2 < 2p(p-1) -2$, we apply Proposition 5.2. If $\pi_{2n + 2p -1}(S^{2m + 1})$ is nontrivial, then there must be an integer $k > m$ such that
$$ 2n - 2m + 2p -2 = 2k(p -1) -2.$$
Thus $(n -m) > m(p-1) -p$, and hence $n > p(m -1)$. Since $p > n$, we have $m = 1$. Therefore, if $m > 1$ then the homotopy group $\pi_{2n + 2p -1}(S^{2m +1})$  is trivial, so we see that $\pi_{2n + 2p -1}(B_m(p))$ is trivial from the fibration (5.3). Thus $[S^{2n + 1} \cup_{\alpha_1} e^{2n + 2p -1}, X] = 0$ for $X = S^{2m + 1}$ or $B_n(p)$  for $m > 1$ by (5.4). 
If $m =1$, then $n - 1 = k(p-1) - p$. Since $p > n$, we have $p - 1 > k(p - 1) - p$, so $k = 2$ and $n = p -1$. Since $\pi_{2p}(S^3)$ is isomorphic to ${\bf Z}_p$ generated by $\alpha_1$, $\pi_{4p-3}(S^3)$ is isomorphic to ${\bf Z}_p$ generated by $\alpha_1^2$ \cite{to}, the homomorphism
$$[S^{2p}, S^3] \xrightarrow{\alpha_1^*} [S^{4p-3}, S^3] $$
is an epimorphism, thus we obtain that $ [S^{2p - 1} \cup_{\alpha_1} e^{4p -3}, S^3] = 0$ by (5.4). Finally, $\pi_{4p -3}(B_1(p))$ is trivial by the fibration (see \cite{oka-1}):  
$$ S^{3} \to B_1(p) \to S^{2p+ 1}.$$
Hence $[S^{2p-1} \cup_{\alpha_1} e^{4p -3}, B_1(p)] = 0$. 
Now Lemma 5.4 follows from Lemma 5.3.    
\end{proof}
The following result follows from the results of \cite{ja2}.
\begin{lemma}[{\rm \cite{mito2}}]
$\Sigma^2B_n(p) \simeq S^{2n + 3} \cup _{\alpha_1} e^{2n + 2p + 1} \vee S^{4n + 2p + 2}.$ 
\end{lemma}
By using Lemma 5.5 we obtain the next lemma which will play an essential role in the remaining of this section.
\begin{lemma}
If $\pi_{4n + 4p -3}(X) = 0$, then ${\rm Im}q^* \subset {\mathcal Z}^\infty(B_n(p), X)$ for a space $X$, where $ q^* : [S^{4n + 2p}, X] \to [B_n(p), X]$.
\end{lemma}
\begin{proof}
Our argument is analogous to that of the proof of Theorem 4.3.
There exists a commutative diagram similar to (4.1).\\
\begin{equation}
\begin{diagram}
\node{B_n(p)} \arrow{e,t}{(1 \times p_{n})\Delta} \arrow{se,r}{((1 \times p_{n})\Delta)'}  
\node{B_n(p) \times S^{2n + 2p -1}} \arrow{e,t}{\wedge} \node{B_n(p) \wedge S^{2n + 2p -1}} \\
\node{} \node{(B_n(p) \times S^{2n + 2p -1})^{4n + 2p}} \arrow{n} \arrow{e,t}{\wedge'} \node{S^{2n +1} \wedge S^{2n + 2p -1},} \arrow{n} \arrow{n}
\end{diagram}
\end{equation}
where $p_n : B_n(p) \to S^{2n + 2p -1}$ is the bundle projection, $((1 \times p_{n})\Delta)'$ is a map induced by $(1 \times p_{n})\Delta$.
If $\pi_{4n + 4p -3}(X) = 0$, any map $ a : S^{4n + 2p} \to X$ extends to $B_n \wedge S^{2n + 2p -1} \simeq S^{4n + 2p} \cup_{\alpha_1} e^{4n + 4p -2} \vee S^{6n + 4p -1}$. The lower map $\wedge'((1 \times p_{n})\Delta)' : B_p(n) \to S^{4n + 2p}$ is homotopic to $\epsilon\cdot q : B_n(p) \to S^{4n + 2p}$, where $\epsilon$ is an integer prime to $p$.
\begin{equation}
\begin{diagram}
\node{B_n(p)} \arrow{e,t}{\wedge(1 \times p_{n})\Delta} \arrow{se,r}{\wedge'((1 \times p_{n})\Delta)'} \node{S^{4n + 2p} \cup_{\alpha_1} e^{4n + 4p -2} \vee S^{6n + 4p -1}} \arrow{se,t}{\overline{a} \vee 0}\\
\node{} \node{S^{2n+1} \wedge S^{2n + 2p -1}} \arrow{n} \arrow{e,t}{a} \node{X}
\end{diagram}
\end{equation}
The upper composition in (5.6) induces the zero on homotopy groups as $\wedge$ does.  Therefore the lower composition in (5.6) induces the trivial map on homotopy groups and so $q^*(a) \in {\mathcal Z}^\infty(B_n(p), X)$.
\end{proof}
\begin{lemma}
Let $p$ be an odd prime, $n$ and $m$ be positive integers with $n \not= m$.
\begin{align*}
(1)~~& sz_p(B_n(p)) = 2n + 2p -1~~~~~\text{and}~~~~lz_p(B_n(p)) = 2~~~~\text{for}~~~~~n \leq 2p - 3,~ p \geq 5,\\
& \text{and for}~~~~n \leq 2, p = 3.\\ 
(2)~~& {\mathcal Z}^\infty(B_n(p), B_m(p)) = [B_n(p), B_m(p)]~~~~~\text{for }~~~~~m, n < p,~ p\geq 3.\\
(3)~~& {\mathcal Z}^\infty(B_n(p), S^{2m +1}) = [B_n(p), S^{2m + 1}]~~~~\text{for } n <p, n-m + p -1 \not= 0,
m \geq 3,\\
& p \geq 5.\\
(4)~~& {\mathcal Z}^\infty(B_1(p), S^{5}) = [B_1(p), S^{5}]~~~~\text{for}~~~~p \geq 3.
\end{align*} 
\end{lemma}
\begin{proof}
(1). We consider the diagram as follows.
\begin{equation}
\begin{diagram}
\node{} \node{} \node{[S^{2n + 1}, B_n(p)]} \\
\node{[S^{4n + 2p}, B_n(p)]} \arrow{e,t}{q^*} \node{[B_n(p), B_n(p)]} \arrow{e,t}{j^*} \node{[S^{2n + 1} \cup_{\alpha_1} e^{2n + 2p -1}, B_n(p)]} \arrow{n} \\
\node{} \node{} \node{[S^{2n + 2p -1}, B_n(p)]} \arrow{n} 
\end{diagram}
\end{equation}
$[S^{2n + 1}, B_n(p)] = \pi_{2n + 1}(B_n(p))$ is isomorphic to ${\bf Z}_{(p)}$, the integers localized at $p$. And $[S^{2n + 2p -1}, B_n(p)] = \pi_{2n + 2p -1}(B_n(p))$ is also isomorphic to ${\bf Z}_{(p)}$. Easily we see that 
$$ [S^{4n + 2p}, B_n(p)] \xrightarrow{q^*} {\mathcal Z}^{2n + 2p -1}(B_n(p))$$
is surjective. If $n \leq 2p -3$ and $p \geq 5$, then $\pi_{4n + 4p - 3}(B_n(p))$ is trivial by Proposition 5.2, and hence we obtain (1) by Lemma 5.6. By (5.7) $lz_p(B_n(p)) = lz_0(B_n(p))$ and the assertion for $lz_0(B_n(p))$ is obtained by Lemma 2.4. Since we see easily that $\pi_{13}(B_1(3)) = 0$ and $\pi_{17}(B_2(3)) = 0 $, the $p = 3$ case is similar.\\ 
(2). According to Oka\cite{oka-1}$, \pi_{4n + 4p -3}(B_m(p))= 0$ for $m, n < p$ and $p \geq 3$. If $n \not= m$, we obtain (2) from Lemma 5.4 and Lemma 5.6.\\
(3). As $4n + 4p -2m -4 <2p(p-1) -2$ for $m \geq 3$, we apply Proposition 5.2 and obtain (3)  by Lemma 5.4 and Lemma 5.6.\\
(4). (4) also follows from Proposition 5.2, Lemma 5.4 and Lemma 5.6.  
\end{proof}
Now we will prove the main theorem in this section.
\begin{proof} (Proof of Theorem 5.1).  By Theorem 4.2\cite{mito2}, $G$ is quasi $p$-regular if and only if
\begin{align*}
 p& > n & for&~~~~Sp(n)\\
 p& >\frac{n}{2}&for&~~~~SU(n)\\
 p& >\frac{n-1}{2}&for&~~~~Spin(n)\\
 p& \geq 5&for&~~~~G_2,~ F_4,~ E_6\\
 p& \geq 11&for&~~~~E_7,~E_8.   
\end{align*}
 The assertion for classical groups $SU(n)$, $Sp(n)$ follows from Lemma 5.7, Lemma 2.4 and the following homotopy equivalences \cite{oka-1}.
$$ \prod_{k=1}^n B_k(p)_p \times \prod_{k=n+1}^{p-1} S^{2k +1}_p \cong SU(n + p)_p\hspace{0.3cm}for~~~~n < p,$$
$$  \prod_{k=1}^n B_{2k-1}(p)_p \times \prod_{k=n+1}^{(p-1)/2} S^{4k - 1}_p \cong Sp(n + (p-1)/2)_p~~~~~for~~~~n < (p-1)/2.$$
  For the $Spin(n)$ case, we use the following equivalences for an odd prime $p$ (see \cite{mito2}).
$$ Spin(2n+1)_p \cong Sp(n)_p$$
$$ Spin(2n)_p \cong Spin(2n -1)_p \times S^{2n-1}_p.$$
 Moreover, our assertion for the exceptional Lie groups follows from Lemma 5.7, Lemma 2.4 and by the results of \cite{mito2} except for the cases $G = F_4, p = 5$, $G =  E_6, p = 5$, $G = E_8, p =11, p = 13$. For these cases there exist equivalences  as follows \cite{mito}.
$$ (F_4)_5 \cong (B_1(5) \times B_7(5))_5 $$
$$ (E_6)_5 \cong (B_1(5) \times B_4(5) \times B_7(5))_5 $$
$$ (E_8)_{11} \cong (B_1(11) \times B_7(11) \times B_{13}(11) \times B_{19}(11))_{11}$$
$$ (E_8)_{13} \cong (B_1(13) \times B_7(13) \times B_{11}(13) \times B_{17}(13))_{13}$$
 We easily see that
all sets $ [B_n(5), B_m(5)],~~~~n \not= m,~~~n, m \in \{1, 4, 7\}$ are trivial except for the case where $n = 4, m = 1$. Thus ${\mathcal Z}^\infty(F_4)_5 = {\mathcal Z}^{23}(F_4)_5$ by Lemma 5.7 (1) and Lemma 2.4. On the other hand, $[B_4(5), B_1(5)]$ is equal to $q^*\pi_{26}(B_1(5))$ which is nontrivial ($\cong {\bf Z}_5$). According to \cite{oka-1}, $\pi_{33}(B_1(5)) = 0$ and we apply Lemma 5.6 and obtain   ${\mathcal Z}^\infty(B_4(5), B_1(5)) = [B_4(5), B_1(5)]$. By Lemma 5.7 (1) and the last equality, we obtain  ${\mathcal Z}^\infty(E_6)_5 = {\mathcal Z}^{23}(E_6)_5$.    
 We have $[B_n(11), B_m(11)] = 0$ for $n, m$ such that $n \not= m$ and $n, m \in \{1, 7, 13, 19 \}$ except for $[B_{13}(11), B_7(11)]$ and also we have  $[B_n(13), B_m(13)] = 0$ for $n, m$ such that $n \not= m$ and $n, m \in \{1, 7, 11, 17\}$ except for $[B_{17}(13), B_{11}(13)]$. The sets $[B_{13}(11), B_7(11)]$ and $[B_{17}(13), B_{11}(13)]$ are both $q^*$-images. Since $\pi_{93}(B_7(11)) = 0 = \pi_{117}(B_{11}(13))$, ${\mathcal Z}^\infty(B_{13}(11), B_7(11)) = [B_{13}(11), B_7(11)]$ and ${\mathcal Z}^\infty(B_{17}(13), B_{11}(13)) = [B_{17}(13), B_{11}(13)]$ by Lemma 5.6, and hence we obtain the $E_8$ case. For the assertion for $lz_p(G)$, we should note that in any quasi p-regular type of our $G$:
$$  G_p \simeq \prod S^{m_i}_p \times \prod B_{n_j}(p)_p, $$  
all homotopy sets $[S^{m_i}, B_{n_j}]$ are trivial for dimensional reasons. Therefore we obtain the result for $lz_p(G)$ by Lemma 2.4 and Lemma 5.7 (1) (and its proof). 
\end{proof}
We close this section with a question on the stability.\\
\ \ \\
{\bf Question.} Is Theorem 5.1 generalized to all $p$ ? More precisely, for a compact, connected simply connected simple Lie group with $H^*(G;{\bf Q}) \cong \Lambda(x_1,...,x_r)$ with ${\rm deg}x_i = n_i$, $n_1 \leq \cdots \leq n_r$, $sz(G) = n_r ?$
\section{Application}
In this section, we will study the groups of self-homotopy equivalences of Lie groups and their subgroups associated to homotopy groups. We will apply the results in the previous sections to determine these groups.  
\begin{definition}
Let us denote by ${\mathcal E}_\#^n(X)$ the group of homotopy classes of self-homotopy equivalences  which induce the identity map on $\pi_i(X)$ for $i \leq n$. 
\end{definition}
\begin{theorem}
${\mathcal E}_\#^\infty(SU(3)) = {\mathcal E}_\#^{n}(SU(3)) \cong {\bf Z}_{12}$ for $n \geq 5$, ${\mathcal E}_\#^\infty(Sp(2)) = {\mathcal E}_\#^n(Sp(2)) \cong {\bf Z}_{120}$ for $n \geq 7$.
\end{theorem}     
\begin{proof}
There exists a bijection
\begin{equation}
T: {\mathcal Z}^n(X) \to {\mathcal E}_\#^n(X) 
\end{equation}
defined by $f \to 1 + f$ for a homotopy associative $H$-space $X$. Thus by Theorem 3.3 ${\mathcal E}_\#^n(SU(3)) = {\mathcal E}_\#^{n+1}(SU(3))$ for $n \geq 5$ and  ${\mathcal E}_\#^n(Sp(2)) = {\mathcal E}_\#^{n+1}(Sp(2))$ for $n \geq 7$. The group structures are obtained by the results of \cite{saw} (or \cite{mio}). 
\end{proof}
More generally we obtain the following theorem which is corresponding to Theorem 5.1. We note that ${\mathcal E}_\#^n(X)$, $n \geq \dim X$, is a nilpotent group for a finite nilpotent space by \cite{dz}.  
\begin{theorem}
Let $G$ be a compact connected, simply-connected simple Lie group,  $H^*(G;{\bf Q}) \cong \Lambda(x_1,...,x_r)$ with ${\rm deg}x_i = n_i$, $n_1 \leq \cdots \leq n_r$. Assume that $G$ is quasi $p$-regular.
\begin{enumerate}
\item[(1)] $ {\mathcal E}_\#^\infty(G)_p = {\mathcal E}_\#^{{\rm dim} G}(G)_p$.\\
\item[(2)] Moreover, if $G = SU(n)$ or $Sp(n)$, then ${\mathcal E}_\#^{n_r}(G)$ is a nilpotent group and ${\mathcal E}_\#^\infty(G)_p = {\mathcal E}_\#^{n_r}(G)_p$.
\end{enumerate}
\end{theorem}
\begin{proof}
(1) If $n \geq {\rm dim}G$, then ${\mathcal E}_\#^{n}(G)_p \cong {\mathcal E}_\#^{n}(G_p)$, see \cite{ma}. Thus the result follows from the bijection $T$ in (6.1) and Theorem 5.1.\\
(2) If $G = SU(n)$ or $Sp(2)$, then $H_*(G)$ is torsion free. Therefore ${\mathcal E}_\#^{n_r}(G)$ acts on homology nilpotently, hence by a result of \cite{dz} it is a nilpotent group. If ${\mathcal E}_\#^{n_r}(G)_p \cong {\mathcal E}_\#^{n_r}(G_p)$, we obtain  
$$ {\mathcal E}_\#^{\infty}(G)_p = {\mathcal E}_\#^{\infty}(G_p) = {\mathcal E}_\#^{n_r}(G_p) = {\mathcal E}_\#^{n_r}(G)_p.$$ 
Thus it suffices to show the following assertion.\\

{\bf Assertion.} The natural map ${\mathcal E}_\#^{n_r}(G) \to {\mathcal E}_\#^{n_r}(G_p)$ is the localization for the above prime $p$.
\\ 
\ \ \\
In the remaining we will prove this assertion. Consider the following commutative diagram:
\begin{equation}
\begin{diagram}
\node{{\mathcal E}_\#^{n_r}(G)} \arrow{e} \node{{\mathcal E}_\#^{n_r}(G_p)} \\
\node{{\mathcal E}_\#^{\infty}(G)} \arrow{e} \arrow{n,l}{j_*} \node{{\mathcal E}_\#^{\infty}(G_p)} \arrow{n,r}{j_*} \arrow{e,t}{\cong} \node{{\mathcal E}_\#^{\dim G}(G_p)} \arrow{nw,r}{\cong}
\end{diagram}
\end{equation}
As we saw in the proof of (1), $j_* : {\mathcal E}_\#^{\infty}(G_p) \to {\mathcal E}_\#^{n_r}(G_p)$ is an isomorphism, and ${\mathcal E}_\#^{\infty}(G) \to {\mathcal E}_\#^{\infty}(G_p)$ is the localization. Thus by (6.2), ${\mathcal E}_\#^{n_r}(G)_p \to {\mathcal E}_\#^{n_r}(G_p)$ is an epimorphism. We also have the following commutative diagram: 
\begin{equation}
\begin{diagram}
\node{{\mathcal Z}^{n_r}(G_p)} \arrow{e,t}{T} \node{{\mathcal E}_\#^{n_r}(G_p)} \\
\node{{\mathcal Z}^{n_r}(G)} \arrow{e,t}{T} \arrow{n} \node{{\mathcal E}_\#^{n_r}(G)} \arrow{n}
\end{diagram}
\end{equation}
where the vertical maps are homomorphisms induced by $p$-localization, the left map is the localization at $p$ by Lemma 4.2 (2). Let $q$ be a prime number and $S_q$ be the $q$-Sylow subgroup of the subgroup of all the torsion elements of ${\mathcal Z}^{n_r}(G)$. We have the following equalities.
$$ T(x)\circ T(y) = (1 + y + x\circ(1 + y)),$$
$$ T((-x)\circ (1 + x)^{-1})\circ T(x) = 1.$$
For any $x, y \in  S_q$. Thus $T(S_q)$ is a subgroup of  ${\mathcal E}_\#^{n_r}(G)$ of order $|S_q|$ since $T$ is bijective. Moreover if $q_1$ and $q_2$ are different primes, then $T(S_{q_1} \times S_{q_2})$ is a subgroup of order $|S_{q_1} \times S_{q_2}|$. Thus, if $x$ is a torsion element of order $|x|$ and $|x|$ is prime to $p$, then the order of $T(x)$ is prime to $p$. It follows that the kernel of the homomorphism  ${\mathcal E}_\#^{n_r}(G) \to {\mathcal E}_\#^{n_r}(G_p)$ is a finite subgroup whose order is prime to $p$, so the map is a localization map at $p$.
\end{proof}
{\bf Remark.} Arkowitz and Strom \cite{as} computed ${\mathcal E}^\infty_\#(G)_p$ when $G$ is $p$-regular in many cases.
%


\begin{thebibliography}{99}
%
\bibitem{as}
M.~Arkowitz and J.~Strom.
\newblock The group of homotopy equivalences of products of spheres and of Lie groups.
\newblock preprint. 
\bibitem{ch}
D. Christensen.
\newblock Ideals in triangulated categories:phantoms, ghosts and skeleta.
\newblock Adv. Math. \textbf{136}, 1998, 284-339.
\bibitem{dz}
E.~Dror and A.~ Zabrodsky 
\newblock  Unipotency and nilpotency in homotopy equivalences. 
\newblock  Topology \textbf{18}, 1979, 187-197. 
\bibitem{fr}
P.~Freyd.
\newblock Stable Homotopy (Proc. Conf. on Cat. Alg.).
\newblock Springer, 1966, 121-172.
\bibitem{ja}
I.~M~.~James.
\newblock On $H$-spaces and their homotopy groups.
\newblock Quart. J. Math. \textbf{11}, 1960, 161-179. 
\bibitem{ja2}
I.~M.~James.
\newblock On sphere bundles over spheres. 
\newblock Comment. Math. Helv. \textbf{35}, 1961, 126-135. 
\bibitem{ku}
P.~G.~Kumpel, Jr.
\newblock Lie groups and products of spheres.
\newblock Proc. Am. Math. Soc.\textbf{16}, 1965, 1350-1356.
\bibitem{ma}
K.~Maruyama.
\newblock Stability properties of maps between Hopf spaces.
\newblock Quart. J.  Math. \textbf{53}, 2002, 47-57.
%
\bibitem{mi2}
M.~Mimura.
\newblock The homotopy groups of Lie groups of low rank.
\newblock  J. Math. Kyoto Univ.  \textbf{6}, 1967, 131-176.
%
\bibitem{mi}
M.~Mimura.
\newblock On the number of multiplication on SU(3) and Sp(2).
\newblock Trans. Amer. Math. Soc.\textbf{146}, 1969, 473-492.
%
\bibitem{mio}
M.~Mimura and H.~\=Oshima.
\newblock Self homotopy groups of {Hopf} spaces with at most three cells.
\newblock J. Math. Soc. Japan \textbf{51}, 1999, 71-92.
\bibitem{mito}
M.~Mimura and H.~Toda.
\newblock Homotopy groups of $SU(3),~SU(4)$ and  $Sp(2)$.
\newblock  J. Math. Kyoto Univ. \textbf{3}, 1964, 217-250.
%
%
\bibitem{mito2}
M.~Mimura and H.~Toda.
\newblock Cohomology operations and the homotopy of compact Lie groups-I.
\newblock  Topology \textbf{9}, 1970, 317-336.
%
%
%
\bibitem{oka-1}
S.~Oka.
\newblock The homotopy groups of sphere bundles over spheres.
\newblock J. Sci. Hiroshima Univ., \textbf{33}, 1969, 161-195.
%
\bibitem{oka}
S.~Oka.
\newblock Homotopy of the exceptional Lie group $G_2$.
\newblock Proc. Edinburgh Math. Soc., \textbf{29}, 1986, 145-169.
%
\bibitem{os}
H.~O\=shima.
\newblock Self homotopy group of the exceptional Lie group $G_2$.
\newblock J. Math. Kyoto Univ \textbf{40}, 2000, 177-184.
%
\bibitem{saw}
N.~Sawashita.
\newblock On $H$-equivalences of $SU(3)$, $U(3)$ and $Sp(2)$.
\newblock J. Math. Tokushima Univ., \textbf{11}, 1977, 33-47.
%
\bibitem{se}
J.~P.~Serre.
\newblock Groupes d'homotopie et classes des groupes ab\'eliens.
\newblock Ann. Math. \textbf{58}, 1953, 258-294.
%
%
\bibitem{to}
H.~Toda.
\newblock Composition methods in homotopy groups of spheres.
\newblock Princeton University Press, Princeton, 1962. 
%
\bibitem{wh}
G.~W.~Whitehead.  
Elements of homotopy theory, 
GTM, 
\textbf{61},  
Springer-Verlag, 
1978 
%
%
\end{thebibliography}
\end{document}